\documentclass[12pt]{article}

\usepackage[francais,english]{babel}
\usepackage{amsfonts,latexsym,amstext}
\usepackage{amsmath,amscd,amssymb}

\usepackage[latin1]{inputenc}

\usepackage{aeguill}

\usepackage[T1]{fontenc}

\usepackage[usenames]{color}

\usepackage{url}

\usepackage{psfrag}

\usepackage{placeins}

\usepackage[dvips]{graphicx}
\usepackage{latexsym,amssymb,amsmath,color,times}
\usepackage{enumerate}

\newtheorem{thm}{Theorem}[section]

\newtheorem{lem}[thm]{Lemma}
\newtheorem{prop}[thm]{Proposition}
\newtheorem{rem}[thm]{Remark}

\numberwithin{equation}{section}
\newcommand{\norm}[1]{\left\Vert#1\right\Vert}
\newcommand{\abs}[1]{\left\vert#1\right\vert}
\newcommand{\set}[1]{\left\{#1\right\}}

\def\R{\mathbb R}
\def\N{\mathbb N}

\def\E{\mathbb E}
\def\P{\mathbb P}
\def\Q{\mathbb Q}

\def\cprime{$'$}

\newcommand{\indicatrice}{\mathchoice{\rm 1\mskip-4mu l}{\rm 1\mskip-4mu l}{\rm 1\mskip-4.5mu l} {\rm 1\mskip-5mu l}}

\def\sqw{\hbox{\rlap{\leavevmode\raise.3ex\hbox{$\sqcap$}}$%
\sqcup$}}
\def\sqb{\hbox{\hskip5pt\vrule width4pt height6pt depth1.5pt%
\hskip1pt}}

\def\qed{\ifmmode\hbox{\hfill\sqb}\else{\ifhmode\unskip\fi%
\nobreak\hfil
\penalty50\hskip1em\null\nobreak\hfil\sqb
\parfillskip=0pt\finalhyphendemerits=0\endgraf}\fi}
\def\cqfd{\ifmmode\sqw\else{\ifhmode\unskip\fi\nobreak\hfil
\penalty50\hskip1em\null\nobreak\hfil\sqw
\parfillskip=0pt\finalhyphendemerits=0\endgraf}\fi}

\begin{document}

\renewcommand{\labelitemi}{$\bullet$}
\bibliographystyle{plain}
\pagestyle{headings}
\title{On the uniqueness of solutions to quadratic BSDEs with convex generators and unbounded terminal conditions: the critical case}

\author{
Freddy Delbaen\\
Department of Mathematics\\
ETH-Zentrum, HG G 54.3, CH-8092 Z\"urich, Switzerland\\
e-mail: delbaen@math.ethz.ch\\ \\
Ying Hu\\
IRMAR, Universit\'e Rennes 1\\ Campus de Beaulieu, F-35042 Rennes
Cedex, France\\ e-mail: ying.hu@univ-rennes1.fr\\ \\
Adrien Richou\\
Univ. Bordeaux, IMB, UMR 5251, F-33400 Talence, France\\
INRIA, \'Equipe ALEA, F-33400 Talence, France\\ 
e-mail: adrien.richou@math.u-bordeaux1.fr}

\date{}
\selectlanguage{english}

\maketitle

\begin{abstract} 

In F.~Delbaen, Y.~Hu and A.~Richou (Ann. Inst. Henri Poincar\'e Probab. Stat. 47(2):559--574, 2011), the authors proved that uniqueness of solution to quadratic BSDE with convex generator and unbounded terminal condition holds  among solutions whose exponentials are $L^p$ with $p$ bigger than a constant $\gamma$ ($p>\gamma$).
In this paper, we consider the critical case: $p=\gamma$. We prove that the uniqueness holds among solutions whose exponentials are $L^\gamma$ under the additional assumption that the generator is
strongly convex. These exponential moments are natural as they are given by the existence theorem.

\end{abstract}

%
%



\section{Introduction}

Since the seminar papers of Bismut \cite{Bismut} and Pardoux-Peng \cite{PP-90}, backward stochastic differential equations (BSDEs in short for the remaining of the paper) have become an active domain of 
research. This is due to, one the one hand, the deep link between BSDEs and Partial Differential Equations (PDEs) (see, e.g., \cite{PP2}), and one the other hand, the profound applications of BSDEs to
mathematical finance and stochastic control theory (see, e.g., the survey paper \cite{EPQ}). A lot of efforts have been made in order to study the well posedness of these equations. We refer to the survey paper \cite{EPQ} and the paper \cite{BDHPS} for existence and uniqueness results about BSDEs with Lipschitz generators.

Quadratic BSDE is a kind of BSDE which has attracted particular attention recently as it appears naturally in utility maximization problems (see, e.g. \cite{RE} and \cite{HIR}) and it is the subject of
the paper. Quadratic BSDEs with bounded terminal conditions are first studied by Kobylanski \cite{Kobylanski-00}. She established rather general results of existence and uniqueness concerning the scalar-valued quadratic BSDEs with bounded terminal conditions. Her proof was simplified by Tevzadze \cite{Tevzadze} as well as by Briand and Elie \cite{BrEl} in some less general cases.
Mania and Schweizer \cite{MS} and  Morlais \cite{Morlais} studied quadratic BSDEs driven by a continuous martingale and applied their results to study some related utility maximization problem. Morlais \cite{Morlais2,Morlais3} considered some quadratic BSDEs with jumps and studied utility maximization in a jump market model.

Concerning the quadratic BSDEs with unbounded terminal value, Briand and Hu \cite{Briand-Hu-06} first showed the existence of solution. Barrieu and El Karoui \cite{BE} revisited the existence issue by a direct forward method that doe not use the result of Kobylanksi.

In this article, we consider the following quadratic BSDE 
\begin{equation}
\label{EDSR}
 Y_t=\xi-\int_t^T g(Z_s)ds +\int_t^T Z_s dW_s, \quad 0 \leqslant t \leqslant T,
\end{equation}
where the generator $g$ is a continuous real function that is convex and has a quadratic growth with respect to the variable $z$. Moreover $\xi$ is an unbounded random variable. Let us recall that, in the previous equation, we are looking for a pair of processes $(Y,Z)$ which is required to be adapted with respect to the filtration generated by the $\mathbb{R}^d$-valued Brownian motion $W$.

In order to state the main result of this paper, let us suppose that there exists a constant $\gamma>0$ such that 
$$\xi^+ \in L^1,  \exp(-\gamma\xi) \in L^1 \mbox{ and } 0\leqslant g(z)\leqslant \frac{\gamma}{2}\abs{z}^2.$$

By a localization procedure similar to that in \cite{Briand-Hu-06}, we prove easily that 
the BSDE (\ref{EDSR}) has at least a solution $(Y,Z)$ such that $e^{-\gamma Y}$ and $Y$ belong to the class (D).

Concerning the uniqueness issue, in \cite{Briand-Hu-08}, the authors proved that the uniqueness holds among solutions whose exponentials are in any $L^p$.
In \cite{DHR-11}, the authors proved that the uniqueness holds among solutions whose exponentials are in $L^p$ for a given $p>\gamma$, i.e. 
$$\mathbb E\left[\sup_{0\le t\le T} e^{pY_t^-}\right]<\infty.$$

However, if we take $g(z)= \frac{\gamma}{2}\abs{z}^2$, then it is easy to see that for the associated BSDE, the uniqueness holds among solutions $(Y,Z)$ such that $e^{-\gamma Y}$ belongs to the class (D). It suffices to note that if $(Y,Z)$ is a solution such that $e^{-\gamma Y}$ belongs to the class (D), then $e^{-\gamma Y}$ is a uniformly integrable martingale and
$$Y_t=-\frac{1}{\gamma} \ln \mathbb{E} \left[\left.e^{-\gamma \xi}\right| \mathcal{F}_t\right].$$

So the aim of this paper is to study the uniqueness of solution of BSDE (\ref{EDSR}) in the critical case: $p=\gamma$. We prove that 
the BSDE (\ref{EDSR}) has a unique solution $(Y,Z)$ such that $e^{-\gamma Y}$ belongs to the class (D) under the additional assumption that the generator $g$ is strongly convex.
We do not know if this result stays true  without this additional assumption.

Let us mention that Richou \cite{Richou-12} and Masiero and Richou \cite{MR}  studied similar problems but they supposed that the BSDEs are Markovian BSDEs.
Cheridito and Nam \cite{CN1} studied quadratic BSDEs with terminal conditions that have bounded Malliavin derivative.
There are also some studies concerning the numerical simulations of quadratic BSDEs (see, e.g., \cite{Richou-11} and \cite{CR}). Finally we note that the works mentioned above are concentrated on scalar-valued quadratic BSDEs.
Concerning multi-dimensional quadratic BSDEs, Tevzadze \cite{Tevzadze} proved that there exists a solution when the terminal value is sufficiently small. Cheridito and Nam \cite{CN2} studied some special Markovian multi-dimensional quadratic BSDEs. There exist some counter-examples in \cite{FR}  showing that multi-dimensional quadratic BSDEs with bounded terminal value may not have a bounded solution.
Kohlmann and Tang \cite{KT}, Tang \cite{Tang}, Hu and Zhou \cite{HZ}, Qian and Zhou \cite{QZ}  proved well-solvability for some special classes of matrix-valued BSDEs.

The paper is organized as follows. Next section is devoted to an existence result, section 3 contains a useful property for solutions and the last section  is devoted to our main uniqueness result.

Let us close this introduction by giving notations that we will use in all the article. For the remaining of the paper, let us fix a nonnegative real number $T>0$. First of all, $(W_t)_{t \in [0,T]}$ is a standard Brownian motion with values in $\mathbb{R}^d$ defined on some complete probability space $(\Omega,\mathcal{F},\mathbb{P})$. $(\mathcal{F}_t)_{t \geqslant 0}$ is the natural filtration of the Brownian motion $W$ augmented by the $\mathbb{P}$-null sets of $\mathcal{F}$. 

As mentioned before, we will deal only with real valued BSDEs which are equations of type (\ref{EDSR}). The function $g$ is called the generator and $\xi$ the terminal condition. Let us recall that a generator is a function $\mathbb{R}^{1\times d} \rightarrow \mathbb{R}$ which is measurable with respect to $ \mathcal{B}(\mathbb{R}^{1\times d})$ and a terminal condition is simply a real $\mathcal{F}_T$-measurable random variable. By a solution to the BSDE (\ref{EDSR}) we mean a pair $(Y_t,Z_t)_{t \in [0,T]}$ of predictable processes with values in $\mathbb{R} \times \mathbb{R}^{1\times d}$ such that $\mathbb{P}$-a.s., $t \mapsto Y_t$ is continuous, $t \mapsto Z_t$ belongs to $L^2(0,T)$, $t \mapsto g(Z_t)$ belongs to $L^1(0,T)$ and $\mathbb{P}$-a.s. $(Y,Z)$ verifies  (\ref{EDSR}).

For any real $p\geqslant 1$, $\mathcal{S}^p$ denotes the set of real-valued, adapted and c\`adl\`ag processes $(Y_t)_{t \in [0,T]}$ such that
$$\norm{Y}_{\mathcal{S}^p}:=\mathbb{E} \left[\sup_{0 \leqslant t\leqslant T} \abs{Y_t}^p \right]^{1/p} < + \infty.$$
$\mathcal{M}^p$ denotes the set of (equivalent class of) predictable processes $(Z_t)_{t \in [0,T]}$ with values in $\mathbb{R}^{1 \times d}$ such that
$$\norm{Z}_{\mathcal{M}^p}:=\mathbb{E}\left[\left(\int_0^T \abs{Z_s}^2 ds \right)^{p/2}\right]^{1/p} < +\infty.$$
We also recall that $Y$ belongs to the class (D) as soon as the family 
$$\set{Y_{\tau}: \tau\leqslant T \textrm{ stopping time}}$$
is uniformly integrable.

For any convex function $f : \R^{1 \times d} \rightarrow \R$, we denote $f^*$ the Legendre-Fenchel transform of $f$ given by
$$f^*(q) = \sup_{z \in \R^{1 \times d}} (zq-f(z)), \quad \forall q \in \mathbb{R}^d.$$ 
We also denote $\partial f$ the subdifferential of $f$. We recall that the subdifferential of $f$ at $z_0$ is the non-empty convex compact set of elements $u \in \R^{d}$ such that
$$f(z)-f(z_0) \geqslant (z-z_0)u, \quad \forall z \in \R^{1 \times d}.$$

Finally, for any predictable process $(q_t)_{t \in [0,T]}$ such that $\int_0^T \abs{q_s}^2ds<+\infty$ $\P$-a.s., we denote $\mathcal{E}(q)$ the Dol\'eans-Dade exponential 
$$\left( \exp \left( \int_0^t q_s dW_s-\frac{1}{2} \int_0^t \abs{q_s}^2 ds \right) \right)_{t \in [0,T]}.$$

\section{An existence result}

Let us begin by giving some assumptions used in this paper.

\paragraph{Assumption A.}
There exists a constant $\gamma>0$ such that
\begin{enumerate}
 \item $\xi^+ \in L^1$ and $\exp(-\gamma\xi) \in L^1$,
 \item $g : \mathbb{R}^d \rightarrow \R^+$ is a convex function that satisfies
 \begin{enumerate}
 \item $g(0)=0$,
 \item there exists a constant $C_1 \geqslant 0$ such that $\forall z \in \R^{1 \times d}$,
 $$g(z)\leqslant C_1+\frac{\gamma}{2}\abs{z}^2.$$
\end{enumerate}
\end{enumerate}
\paragraph{Assumption B.}
There exist two constants $\varepsilon>0$ and $C_2 \geqslant 0$ such that $\forall z,z' \in \R^{1\times d}$, $\forall s \in \partial g(z')$,
 $$g(z)-g(z')-(z-z')s  \geqslant \frac{\varepsilon}{2} \abs{z-z'}^2-C_2.$$

\begin{rem}
$ $
 \begin{itemize}
  \item If $g$ is a $C^2$ function then assumption B is equivalent to the assumption: there exist $R \geqslant 0$ and $\varepsilon >0$ such that for all $z \in \R^{1 \times d}$ with $\abs{z} >R$, we have  $g''(z) \geqslant \varepsilon Id$.
  \item For a general convex generator $g$ with quadratic growth it is easy to modify the terminal condition and the probability to obtain a new generator $\tilde{g} : \mathbb{R}^{1\times d} \rightarrow \R^+$ such that assumption A.2. holds true.
 \end{itemize}
\end{rem}

%
%

The aim of this section is to show the existence of solutions under the assumption A, using a localization method.

\begin{thm}
\label{thm existence}
 Let us assume that assumption A holds. Then the BSDE (\ref{EDSR}) has at least a solution $(Y,Z)$ such that:
 $$-\frac{1}{\gamma} \ln \mathbb{E} \left[\left.e^{\gamma C_1 T}e^{\gamma \xi^-}\right| \mathcal{F}_t\right]  \leqslant Y_t \leqslant \E \left[\left.\xi\right| \mathcal{F}_t\right].$$ 
 In particular, $e^{-\gamma Y}$ and $Y$ belong to the class (D).
\end{thm}
\paragraph{Proof of Theorem \ref{thm existence}.}
To show this existence result we use the same classical localization argument as Briand and Hu in \cite{Briand-Hu-06}. Let us fix $n,p \in \N^*$ and set $\xi^{n,p}=\xi^+ \wedge n -\xi^- \wedge p$. Then it is known from \cite{Kobylanski-00} that the BSDE
$$Y_t^{n,p} = \xi^{n,p} -\int_t^T g(Z^{n,p}_s)ds +\int_t^T Z^{n,p}_s dW_s, \quad 0 \leqslant t \leqslant T,$$
has a unique solution $(Y^{n,p},Z^{n,p}) \in \mathcal{S}^{\infty} \times \mathcal{M}^2$. By applying Theorem $2$ in \cite{Briand-Hu-06}, we have the estimate
$$-\frac{1}{\gamma} \ln \mathbb{E} \left[\left. \phi_t(-\xi^{n,p})\right| \mathcal{F}_t\right]  \leqslant Y_t^{n,p}$$
where $(\phi_t(z))_{t \in [0,T]}$ stands for the solution to the integral equation
$$\phi_t(z)=e^{\gamma z} + \int_t^T H(\phi_s(z))ds, \quad 0 \leqslant t \leqslant T,$$
with
$$H(p)=C_1\gamma p \indicatrice_{[1,+\infty[} (p) +C_1 \gamma \indicatrice_{]-\infty,1[} (p).$$
It is noticed in \cite{Briand-Hu-06} that $\phi_t(z) = e^{\gamma C_1 (T-t)}e^{\gamma z}$ when $z \geqslant 0$ and $z \mapsto \phi_t(z)$ is an increasing continuous function. Thus, we have
\begin{eqnarray*}
-\frac{1}{\gamma} \ln \mathbb{E} \left[\left.e^{\gamma C_1 T}e^{\gamma \xi^-}\right| \mathcal{F}_t\right] &\leqslant& -\frac{1}{\gamma} \ln \mathbb{E} \left[\left.e^{\gamma C_1 T}e^{\gamma (\xi^{n,p})^-}\right| \mathcal{F}_t\right]\\
 &\leqslant& -\frac{1}{\gamma} \ln \mathbb{E} \left[\left. \phi_t(-\xi^{n,p})\right| \mathcal{F}_t\right] \leqslant Y_t^{n,p}.
\end{eqnarray*}
Moreover, $g$ is a nonnegative function, so
$$Y_t^{n,p} = \E \left[\left.\xi^{n,p} - \int_t^T g(Z_s^{n,p})ds\right| \mathcal{F}_t\right] \leqslant \E \left[\left.\xi^{n,p} \right| \mathcal{F}_t\right]\leqslant \E \left[\left.\xi^+\right| \mathcal{F}_t\right] .$$ 
We remark that 
$$\forall t \in [0,T], \quad Y_t^{n,p+1} \leqslant Y_t^{n,p} \leqslant Y_t^{n+1,p},$$
and we define $Y^p=\sup_{n \geqslant 1} Y^{n,p}$ so that $Y_t^{p+1} \leqslant Y_t^p$ and $Y_t = \inf_{p \geqslant 1} Y_t^p$. By the dominated convergence theorem, we have
 $$-\frac{1}{\gamma} \ln \mathbb{E} \left[\left.e^{\gamma C_1 T}e^{\gamma \xi^-}\right| \mathcal{F}_t\right]  \leqslant -\frac{1}{\gamma} \ln \mathbb{E} \left[\left. \phi_t(-\xi)\right| \mathcal{F}_t\right]  \leqslant Y_t \leqslant \E \left[\left.\xi\right| \mathcal{F}_t\right],$$
 and in particular, we remark that $\lim_{t \rightarrow +\infty} Y_t = \xi=Y_T$. Arguing as in \cite{Briand-Hu-06} with a localization argument, we can show that there exists a process $Z$ such that $(Y,Z)$ solves the BSDE (\ref{EDSR}). Finally, since processes $t \mapsto \mathbb{E} \left[\left.e^{\gamma C_1 T}e^{\gamma \xi^-}\right| \mathcal{F}_t\right]$ and $t \mapsto \E \left[\left.\xi\right| \mathcal{F}_t\right]$ belong to the class (D), we conclude that $e^{-\gamma Y}$, $Y^+$ and so $Y$ belong to the class (D).

\cqfd

\section{A uniform integrability property for solutions}
In this part we will show the following proposition.
\begin{prop}
\label{expo uniformement integrable}
 We assume that assumption A holds true. Let us consider $(Y,Z)$ a solution of the BSDE (\ref{EDSR}) such that $Y$ and $e^{-\gamma Y}$ belong to the class (D). Then, for all predictable process $(q_s)_{s \in [0,T]}$ with values in $\R^d$ and such that $q_s \in \partial g(Z_s)$ for all $s \in [0,T]$, $\mathcal{E}(q)$ is a uniformly integrable process and defines a probability $\mathbb{Q} \sim \mathbb{P}$.
\end{prop}

\paragraph{Proof of Proposition \ref{expo uniformement integrable}.}
Let us start the proof by giving a simple lemma.
\begin{lem}
\label{delavallePoussin}
 The family of random variables $\set{e^{\gamma X}|X \in \mathcal{H}}$ is uniformly integrable if and only if there exists a function $k : \R^+ \rightarrow \R^+$ such that $k(x) \rightarrow +\infty$ when $x \rightarrow +\infty$, and 
 $$\sup_{X \in \mathcal{H}} \E[K(X^+)] < +\infty,$$
with $K(x)=\int_0^x k(t)e^{\gamma t}dt$. Moreover, we can assume without restriction that $k \in C^{\infty}$, $k(0)=\gamma$ and $k'(x)>0$ for all $x \in \R^+$.
\end{lem}
\paragraph{Proof of Lemma \ref{delavallePoussin}.}
We only prove the nontrivial implication. Firstly, let us remark that $\set{e^{\gamma X}|X \in \mathcal{H}}$ is uniformly integrable if and only if $\set{e^{\gamma X^+}|X \in \mathcal{H}}$ is also uniformly integrable, so we can assume that $\mathcal{H}$ is a family of positive random variables. Now we apply the de la Vall\'ee-Poussin theorem: there exists a nondecreasing function $g :\R^+ \rightarrow \R^+$ which is a constant function on each interval $[n, n+1[$ for $n \in \N$, that satisfies $g(x)\rightarrow +\infty$ when $x \rightarrow +\infty$ and such that
 $$\sup_{X \in \mathcal{H}} \E[G(e^{\gamma X})] < +\infty,$$
 with $G(x)=\int_1^x g(t)dt$. Then, it is simple to consider a smooth approximation $\tilde{g}$ of $g$ such that $\tilde{g}(1)=1$, $\tilde{g}'(x)>0$ for all $x \in [1,+\infty[$ and 
$$g+1-g(1) \leqslant \tilde{g} \leqslant g+C.$$
 This function $\tilde{g}$ also satisfies $\tilde{g}(x)\rightarrow +\infty$ when $x \rightarrow +\infty$ and
  $$\sup_{X \in \mathcal{H}} \E[\tilde{G}(e^{\gamma X})] < +\infty,$$
  with $\tilde{G}(x)=\int_1^x \tilde{g}(t)dt$. A simple calculus gives us
  $$\tilde{G}(e^{\gamma x})=\int_0^t \tilde{g}(e^{\gamma u})\gamma e^{\gamma u} du$$
  and so we just have to set $k(x)= \gamma \tilde{g}(e^{\gamma x})$ to conclude the proof.
  \cqfd

  Now, let us apply the previous lemma in our situation: since we consider a solution $(Y,Z)$ such that $e^{-\gamma Y}$ belongs to the class (D), then there exists a function $k : \R^+ \rightarrow \R^+$ given by Lemma \ref{delavallePoussin} such that
  \begin{equation}
  \label{majoration uniforme EK}
  \sup_{0 \leqslant \tau \leqslant T, \textrm{ stopping time}} \E[K(Y^-_{\tau})] < +\infty,
\end{equation}
  with $K(x)=\int_0^x k(t)e^{\gamma t}dt$. We define
$$\Psi_0(x)=e^{\gamma x}-\gamma x -1=\int_0^x \gamma (e^{\gamma u} -1)du$$ 
and
$$\Psi(x)=\int_0^x k(u)(e^{\gamma u}-1)du.$$
Since $\Psi_0$ and $\Psi$ are convex functions we can also consider their dual functions. $\Phi_0(x)=\left(\frac{x}{\gamma}+1\right) \ln \left( \frac{x}{\gamma}+1 \right) -\frac{x}{\gamma}$ is the dual function of $\Psi_0$ since $\Phi_0'(x)=\frac{1}{\gamma} \ln \left( \frac{x}{\gamma}+1 \right)$ is the inverse function of $\Psi_0'$. Moreover, the dual function of $\Psi$ is given by $\Phi(x)=\int_0^x \Phi'(u)du$ with $\Phi'$ the inverse function of $\Psi'$.

Now we consider a predictable process $(q_s)_{s \in [0,T]}$ with values in $\R^d$ and such that $q_s \in \partial g(Z_s)$ for all $s \in [0,T]$. Firstly let us show that $s \mapsto q_s$ belongs to $L^2(0,T)$ $\P$-a.s.. Since assumption A.2 holds true for $g$, then $g^*$ satisfies
\begin{equation}
 \label{croissance g^*}
 g^*(q) \geqslant -C+\frac{1}{2\gamma} \abs{q}^2 \quad \quad \textrm{and} \quad \quad g^*(0)=0,
\end{equation}
and thus,
\begin{eqnarray*}
\int_0^T \abs{q_s}^2ds \leqslant C+C\int_0^T g^*(q_s)ds& =& C+C\int_0^T (Z_s q_s-g(Z_s))ds\\ &\leqslant& C+C\int_0^T \abs{Z_s q_s}ds+ C\int_0^T \abs{Z_s}^2ds.
\end{eqnarray*}
Moreover, since $q_s \in \partial g(Z_s)$ we have
$$Z_s q_s =(2Z_s -Z_s)q_s \leqslant g(2Z_s)-g(Z_s)$$
and
$$-Z_s q_s =(0 -Z_s)q_s \leqslant g(0)-g(Z_s).$$
So we finally obtain
$$\int_0^T \abs{q_s}^2ds \leqslant  C+ C\int_0^T \abs{Z_s}^2ds<+\infty \quad \P\textrm{-a.s.}.$$

Now let us show that $\mathcal{E}(q)$ is a uniformly integrable martingale. We start by defining the stopping time
$$\tau_n=\inf \set{ t \in [0,T] : \sup\left( \int_0^t \abs{q_s}^2 ds, \int_0^t \abs{Z_s}^2ds \right) \geqslant n} \wedge T,$$
and the probability
$$\frac{d\Q_n}{d\P} = M_{\tau_n}, \quad \textrm{ with } \quad  M_t=\exp\left( \int_0^t q_s dW_s-\frac{1}{2}\int_0^t \abs{q_s}^2 ds\right).$$
We will show that $(M_{\tau_n})_{n \in \mathbb{N}}$ is uniformly integrable which is sufficient to conclude. Since $(Y,Z)$ solves the BSDE (\ref{EDSR}), we have
\begin{eqnarray}
\nonumber
 Y_0& =& Y_{\tau_n} -\int_0^{\tau_n} g(Z_s)ds + \int_0^{\tau_n} Z_s dW_s\\
 \nonumber
 &=& Y_{\tau_n}+ \int_0^{\tau_n} (Z_s q_s -g(Z_s))ds + \int_0^{\tau_n} Z_s (dW_s- q_sds)\\
 \label{Y0=EQ...}
 &=& \E^{\Q_n} \left[Y_{\tau_n} + \int_0^{\tau_n} g^*(q_s)ds\right].
\end{eqnarray}
Firstly, since $\Psi$ and $\Phi$ are dual functions, the Fenchel's inequality gives us
$$ \E^{\Q_n} \left[Y_{\tau_n} \right] \geqslant  -\E^{\Q_n} \left[Y_{\tau_n}^- \right] \geqslant  -\E \left[\Psi(Y_{\tau_n}^-)\right]- \E \left[\Phi(M_{\tau_n})\right].$$
Moreover, we have, thanks to (\ref{majoration uniforme EK}),
$$-\E \left[\Psi(Y_{\tau_n}^-)\right] \geqslant -\E \left[K(Y_{\tau_n}^-)\right]  \geqslant -C,$$
with $C$ a constant that does not depend on $n$. By putting these inequalities into (\ref{Y0=EQ...}) we obtain
\begin{equation}
\label{minorationY0}
Y_0 \geqslant -C - \E \left[\Phi(M_{\tau_n})\right] +\E^{\Q_n} \left[ \int_0^{\tau_n} g^*(q_s)ds\right].
\end{equation}
Thanks to the growth of $g^*$ given by (\ref{croissance g^*}) we have
$$\E^{\Q_n} \left[ \int_0^{\tau_n} g^*(q_s)ds\right] \geqslant -C+\E^{\Q_n} \left[\frac{1}{2\gamma} \int_0^{\tau_n} \abs{q_s}^2ds\right].$$
Moreover, a simple calculus gives us
$$\E^{\Q_n} \left[\frac{1}{2\gamma} \int_0^{\tau_n} \abs{q_s}^2ds\right] = \frac{1}{\gamma} \E \left[M_{\tau_n} \ln  (M_{\tau_n})\right].$$
By putting these two results into (\ref{minorationY0}), and by setting $\Lambda=\Phi_0-\Phi$, we obtain
\begin{equation}
 \label{minorationY02}
 Y_0 \geqslant -C +\E \left[\Lambda(M_{\tau_n})\right]- \E \left[\Phi_0(M_{\tau_n})\right] +\frac{1}{\gamma}\E \left[M_{\tau_n} \ln (M_{\tau_n})\right].
\end{equation}

Let us remark that
\begin{eqnarray*}
 & &\E \left[\Phi_0(M_{\tau_n})\right]-\frac{1}{\gamma}\E \left[M_{\tau_n} \ln (M_{\tau_n})\right]\\
 &=& \E \left[\frac{M_{\tau_n}}{\gamma} \ln \left(1+\frac{\gamma}{M_{\tau_n}}\right)\right]-\left(\frac{\ln \gamma+1}{\gamma}\right)\E \left[M_{\tau_n}\right] + \E \left[ \ln \left(\frac{M_{\tau_n}}{\gamma}+1\right)\right]\\
 &=& \E \left[\frac{M_{\tau_n}}{\gamma} \ln \left(1+\frac{\gamma}{M_{\tau_n}}\right)\right]-\left(\frac{\ln \gamma+1}{\gamma}\right) + \E \left[ \ln \left(\frac{M_{\tau_n}}{\gamma}+1\right)\right].
\end{eqnarray*}
An elementary inequality gives us
$$\E \left[\frac{M_{\tau_n}}{\gamma} \ln \left(1+\frac{\gamma}{M_{\tau_n}}\right)\right] \leqslant \E \left[\frac{M_{\tau_n}}{\gamma} \frac{\gamma}{M_{\tau_n}}\right] \leqslant 1,$$
and
$$\E \left[ \ln \left(\frac{M_{\tau_n}}{\gamma}+1\right)\right] \leqslant \E \left[\frac{M_{\tau_n}}{\gamma}\right] \leqslant \frac{1}{\gamma}.$$
Thus, we have
$$\E \left[\Phi_0(M_{\tau_n})\right]-\frac{1}{\gamma}\E \left[M_{\tau_n} \ln (M_{\tau_n})\right] \leqslant C$$
and inequality (\ref{minorationY02}) becomes
\begin{equation}
 \label{majoration ELambda}
 Y_0 \geqslant -C +\E \left[\Lambda(M_{\tau_n})\right].
\end{equation}
Let us give a useful property of $\Lambda$ that we will prove after.
\begin{prop}
\label{croissance Lambda}
 The function $\Lambda$ satisfies
 $$\lim_{x \rightarrow +\infty} \frac{\Lambda(x)}{x} =+\infty.$$ 
\end{prop}
Thanks to this proposition and the inequality (\ref{majoration ELambda}) we are allowed to apply the de la Vall\'ee-Poussin Theorem: $(M_{\tau_n})_{n \in \mathbb{N}}$ is uniformly integrable and the proof is finished.
\cqfd
\paragraph{Proof of Proposition \ref{croissance Lambda}:}
It is sufficient to show that $\Lambda'=\Phi_0'-\Phi'$ is increasing and $\lim_{x\rightarrow + \infty} \Lambda'(x)\rightarrow +\infty$. 
Firstly, let us show that $\Psi''(\Phi'(x)) \geqslant \gamma( x+\gamma)$, for all $x \geqslant 0$:
\begin{equation*}
 \Psi''(x) = k'(x)(e^{\gamma x}-1)+k(x)\gamma e^{\gamma x} \geqslant \gamma k(x)(e^{\gamma x}-1)+\gamma k(x) \geqslant \gamma \Psi'(x)+\gamma^2,
\end{equation*}
so we have
$$\Psi''(\Phi'(x)) \geqslant \gamma\Psi'(\Phi'(x))+\gamma^2 =\gamma (x+\gamma).$$
As a result, we get from the equality $(\Psi'(\Phi'(x)))'=\Psi''(\Phi'(x))\Phi''(x)=1$ that 
$$\Phi''(x) \leqslant \frac{1}{\gamma(x+\gamma)}.$$ 
We finally obtain
$$\Lambda''(x)=\Phi_0''(x)-\Phi''(x) \geqslant \frac{1}{\gamma(x+\gamma)}-\frac{1}{\gamma(x+\gamma)} \geqslant 0$$
and so $\Lambda'$ is an increasing function. 

To conclude we will prove by contradiction that $\Lambda'$ is an unbounded function: let us assume that there exists a constant $A$ such that $\Lambda' \leqslant A$. Then we have
\begin{eqnarray*}
 x&=& \Psi'(\Phi'(x))=k(\Phi'(x))\left(e^{\gamma\Phi'(x)}-1 \right) = k(\Phi'(x))\left(e^{\gamma(\Phi_0'(x)-\Lambda'(x))}-1 \right)\\
 &\geqslant&  k(\Phi'(x))\left(e^{\gamma\Phi_0'(x)}e^{-A}-1 \right) \geqslant k(\Phi'(x))\left(\left(\frac{x}{\gamma}+1\right)e^{-A}-1 \right),
\end{eqnarray*}
and so, we get for $x$ big enough
$$k(\Phi'(x)) \leqslant \frac{x}{\left(\frac{x}{\gamma}+1\right)e^{-A}-1} \leqslant C.$$
Since $\lim_{x \rightarrow +\infty}\Phi'(x)=+\infty$, previous inequality gives us that $k$ is a bounded function, which is a contradiction.
\cqfd

\begin{rem}\label{entropy}
 In \cite{DHR-11}, the authors proved that if for some $p>\gamma$,
$$\mathbb E\left[\sup_{0\le t\le T} e^{pY_t^-}\right]<\infty,$$
then $\mathcal{E}(q)$ has finite entropy, i.e., 
$$\E \left[ \mathcal{E}(q)_T \ln \mathcal{E}(q)_T\right] <+\infty.$$

However, in the critical case, $e^{-\gamma Y}$ belongs to the class (D), this property is not always true.
It suffices to take again $g(z)= \frac{\gamma}{2}\abs{z}^2$, then if $Y\le 0$,
$$\mathcal{E}(q)_t \ln \mathcal{E}(q)_t= e^{\gamma Y_0}e^{\gamma Y_t^-}(\gamma Y_0+\gamma Y_t^-).$$
It follows that if $g(z)= \frac{\gamma}{2}\abs{z}^2$ and $Y\le 0$, $\mathcal{E}(q)$ has finite entropy if and only if $Y^- e^{\gamma Y^-}$ belongs to the class (D).
\end{rem}

\section{The uniqueness result}

Remark \ref{entropy} indicates that $\mathcal{E}(q)$ does not always have finite entropy in the critical case. Hence we could not adopt the verification argument given in \cite{DHR-11} to show the uniqueness.
In this last section, we show the uniqueness under the additional assumption B.

\begin{thm}
\label{theorem unicite}
 Let us assume that assumptions A and B hold true. Then the BSDE (\ref{EDSR}) has a unique solution $(Y,Z)$ such that $Y$ and $e^{-\gamma Y}$ belong to the class (D).
\end{thm}
\paragraph{Proof of Theorem \ref{theorem unicite}.}
The existence result is already given in Theorem \ref{thm existence}. For the uniqueness, let us consider $(Y,Z)$ and $(Y',Z')$ two solutions of the BSDE (\ref{EDSR}) such that $Y$, $Y'$, $e^{-\gamma Y}$ and $e^{-\gamma Y'}$ belong to the class (D). By a symmetry argument it is sufficient to show that $Y_t \geqslant Y_t'$ $\P$-a.s. for all $ t \in [0,T]$. For $t \in [0,T[$, let us denote $A:= \set{Y_t < Y_t'}$ and set the stopping time $\tau =\inf \set{s \geqslant t| Y_s \geqslant Y'_s}$. Then, for $s \in [t,\tau]$ we have $Y_s \leqslant Y'_s$ and $Y_{\tau}  = Y'_{\tau}$ $\P$-a.s. because $t \rightarrow Y_t$ is continuous $\P$-a.s..

Let us consider a predictable process $(q_s)_{s \in [0,T]}$ with values in $\R^d$ and such that $q_s \in \partial g(Z_s)$ for all $s \in [0,T]$. Thanks to Proposition \ref{expo uniformement integrable} we know that $\mathcal{E}(q)$ defines a probability that we will denote $\mathbb{Q}$. Under $\Q$, we get
\begin{equation}
\label{diff Y-U}
 d(Y_s-Y'_s)=(g(Z_s)-g(Z'_s)-(Z_s-Z'_s)q_s)ds-(Z_s-Z'_s)dW_s^{\Q}.
\end{equation}
Then, Itô formula gives us, for $0<\alpha\leqslant \varepsilon$,
\begin{eqnarray*}
& & de^{\alpha(Y_{s \wedge \tau}-Y_{s \wedge \tau}'+C_2(T-s))\indicatrice_{A}}\\
&=&-\alpha \indicatrice_{A}e^{\alpha((Y_{s \wedge \tau}-Y_{s \wedge \tau}'+C_2(T-s))\indicatrice_{A}}\Big(g(Z'_s)-g(Z_s)-(Z'_s-Z_s)q_s+C_2\\
& &\quad\quad -\frac{\alpha}{2}\abs{Z'_s-Z_s}^2\Big)\indicatrice_{s \leqslant \tau}ds\\
& & -\alpha \indicatrice_{A}e^{\alpha C_2(T-s)\indicatrice_{A}}\indicatrice_{s > \tau}ds +\indicatrice_{s \leqslant \tau}dM_s,
\end{eqnarray*}
with $(M_s)_{s \in [t,\tau]}$ a local martingale under $\mathbb{Q}$. From assumption B we have that
$$g(Z'_s)-g(Z_s)-(Z'_s-Z_s)q_s \geqslant \frac{\varepsilon}{2} \abs{Z'_s-Z_s}^2-C_2.$$
So, we obtain that $\left(e^{\alpha(Y_{s \wedge \tau}-Y_{s \wedge \tau}'+C_2(T-s ))\indicatrice_{A}}\right)_{t \leqslant s \leqslant T}$ is a bounded supermartingale under $\Q$ and
$$e^{\alpha(Y_{s \wedge \tau}-Y_{s \wedge \tau}'+C_2(T-s))\indicatrice_{A}} \geqslant \E\left[\left.e^{\alpha(Y_{\tau}-Y_{\tau}'+C_2(T-T))\indicatrice_{A}}\right| \mathcal{F}_s \right]=1, \quad \forall s \in [t,T].$$
It implies that $((Y_{s \wedge \tau}-Y_{s \wedge \tau}')\indicatrice_{A})_{ s \in [t,T]}$ is a bounded process. Moreover, $g$ is a convex function so
$$g(z)-g(z')-(z-z')u \leqslant 0, \quad \forall z,z' \in \mathbb{R}^{1 \times d}, \quad \forall u \in \partial g(z).$$
By using this inequality in (\ref{diff Y-U}), we obtain that $((Y_{s \wedge \tau}-Y_{s \wedge \tau}')\indicatrice_{A})_{ s \in [t,T]}$ is a bounded negative supermartingale under $\Q$ such that $(Y_{\tau}-Y'_{\tau})\indicatrice_{A}=0$. We conclude that $(Y_{t}-Y'_{t})\indicatrice_{A}=0$, that is to say, $Y_t \geqslant Y'_t$. Finally, it is rather standard to show that $\int_0^T \abs{Z_s-Z_s'}^2ds =0$ $\P$-a.s..
\cqfd


\begin{thebibliography}{}

\end{thebibliography}


\begin{thebibliography}{99}

\bibitem{BE}
P.~Barrieu and N.~El Karoui. 
\newblock Monotone stability of quadratic semimartingales with applications to unbounded general quadratic BSDEs.
\newblock {\em Ann. Probab.}, {\bf 41}(3B):1831--1863, 2013.

\bibitem{Bismut} 
J. M.~Bismut.
\newblock Conjugate convex functions in optimal stochastic control.
\newblock {\em J. Math. Anal. Appl.},  {\bf 44}:384--404, 1973.

\bibitem{BDHPS}
P.~Briand, B.~Delyon, Y.~Hu, E.~Pardoux and L.~Stoica.
\newblock $L^p$ solutions of backward stochastic differential equations.
\newblock {\em Stochastic Process. Appl.}, {\bf 108}(1):109--129, 2003.

\bibitem{BrEl} 
P.~Briand and R. Elie.
\newblock A simple constructive approach to quadratic BSDEs with or without delay.
\newblock {\em Stochastic Process. Appl.}, {\bf 123}(8):2921--2939, 2013.

\bibitem{Briand-Hu-06}
P.~Briand and Y.~Hu.
\newblock B{SDE} with quadratic growth and unbounded terminal value.
\newblock {\em Probab. Theory Related Fields}, {\bf 136}(4):604--618, 2006.

\bibitem{Briand-Hu-08} P.~Briand and Y.~Hu.
\newblock  Quadratic BSDEs with convex generators and unbounded terminal conditions. 
\newblock {\em Probab. Theory Related Fields}, {\bf 141}(3--4):543--567, 2008.

\bibitem{CN1} P.~Cheridito and K.~Nam.
\newblock BSDEs with terminal conditions that have bounded Malliavin derivative.
\newblock arXiv:1211.1089, 2012.


\bibitem{CN2} P.~Cheridito and K.~Nam.
\newblock Multidimensional quadratic and subquadratic BSDEs with special structure.
\newblock arXiv:1309.6716, 2013.

\bibitem{CR} J. F.~Chassagneux and A.~Richou.
\newblock Numerical simulation of quadratic BSDEs.
\newblock arXiv:1307.5741, 2013.


\bibitem{DHR-11} F.~Delbaen, Y.~Hu and A.~Richou.
\newblock  On the uniqueness of solutions to quadratic BSDEs with convex generators and unbounded terminal conditions.
\newblock {\em Ann. Inst. Henri Poincar\'e Probab. Stat.}, {\bf 47}(2):559--574, 2011.

\bibitem{EPQ} N.~El Karoui, S.~Peng and M. C. Quenez.
\newblock Backward stochastic differential equations in finance. 
\newblock {\em Math. Finance}, {\bf 7}(1):1--71, 1997.

\bibitem{FR} C.~Frei and G. dos Reis.
\newblock A financial market with interacting investors: does an equilibrium exist?
\newblock {\rm Math. Financ. Econ.}, {\bf 4}(3):161--182, 2011.

\bibitem{HIR} Y.~Hu, P.~Imkeller and M. Muller.
\newblock Utility maximization in incomplete markets. 
\newblock {\em Ann. Appl. Probab.}, {\bf 15}(3):1691--1712, 2005.

\bibitem{HZ} Y.~Hu and X. Y.~Zhou.
\newblock Indefinite stochastic Riccati equations.
\newblock {\em SIAM J. Control Optim.}, {\bf 42}(1):123--137, 2003.


\bibitem{Kobylanski-00}
M.~Kobylanski.
\newblock Backward stochastic differential equations and partial differential
  equations with quadratic growth.
\newblock {\em Ann. Probab.}, {\bf 28}(2):558--602, 2000.

\bibitem{KT} M.~Kohlmann and S. Tang.
\newblock Multidimensional backward stochastic Riccati equations and applications.
\newblock {\em SIAM J. Control Optim.}, {\bf 41}(6):1696--1721, 2003.


\bibitem{MR} F.~Masiero and A.~Richou.
\newblock A note on the existence of solutions to Markovian superquadratic BSDEs with an unbounded terminal condition.
\newblock {\em Electron. J. Probab.}, {\bf 18}, no. 9, 23 pp., 2013.

\bibitem{MS} M.~Mania and M.~Schweizer.
\newblock Dynamic exponential utility indifference valuation. 
\newblock {\em Ann. Appl. Probab.}, {\bf 15}(3):2113--2143, 2005.

\bibitem{Morlais} M. A. ~Morlais. 
\newblock Quadratic BSDEs driven by a continuous martingale and applications to the utility maximization problem. 
\newblock {\em Finance Stoch.}, {\bf 13}(1):121--150, 2009.

\bibitem{Morlais2} M. A. ~Morlais. 
\newblock Utility maximization in a jump market model.
\newblock {\em Stochastics}, {\bf 81}(1):1--27, 2009.

\bibitem{Morlais3} M. A. ~Morlais. 
\newblock A new existence result for quadratic BSDEs with jumps with application to the utility maximization problem.
\newblock {\em Stochastic Process. Appl.}, {\bf 120}(10):1966--1995, 2010.

\bibitem{PP-90}
E.~Pardoux and S.~Peng.
\newblock Adapted solution of a backward stochastic differential equation.
\newblock {\em Systems Control Lett.}, {\bf 14}(1):55--61, 1990.

\bibitem{PP2} E.~Pardoux and S.~Peng.
\newblock Backward stochastic differential equations and quasilinear parabolic partial differential equations. 
\newblock {\em Stochastic partial differential equations and their applications (Charlotte, NC, 1991), Lecture Notes in Control and Inform. Sci.}, {\bf 176}:200--217, Springer, Berlin, 1992. 

\bibitem{QZ} Z.~Qian and X. Y.~Zhou. 
\newblock Existence of solutions to a class of indefinite stochastic Riccati equations.
\newblock {\em SIAM J. Control Optim.}, {\bf 51}(1):221--229, 2013.

\bibitem{Richou-11} A.~Richou.
\newblock  Numerical simulation of BSDEs with drivers of quadratic growth.
\newblock {\em Ann. Appl. Probab.},  {\bf 21}(5):1933--1964, 2011.

\bibitem{Richou-12} A.~Richou.
\newblock  Markovian quadratic and superquadratic BSDEs with an unbounded terminal condition.
\newblock {\em Stochastic Process. Appl.},  {\bf 122}(9), 3173--3208, 2012.

\bibitem{RE} R.~Rouge and N.~El Karoui.
\newblock  Pricing via utility maximization and entropy. 
\newblock {\rm Math. Finance}, {\bf 10}(2):259--276, 2000.

\bibitem{Tang} S.~Tang.
\newblock General linear quadratic optimal stochastic control problems with random coefficients: linear stochastic Hamilton systems and backward stochastic Riccati equations. 
\newblock {\em SIAM J. Control Optim.}, {\bf 42}(1):53--75, 2003.

\bibitem{Tevzadze} R.~Tevzadze.
\newblock Solvability of backward stochastic differential equations with quadratic growth.
\newblock {\em Stochastic Process. Appl.},  {\bf 118}(3):503--515, 2008.








\end{thebibliography}
\def\cprime{$'$}

\end{document}